\documentclass[11pt]{amsart} 
\usepackage{amsthm,amsbsy,amsmath,amssymb,amscd,amsfonts,array,mathrsfs,verbatim,enumerate,xypic}
\usepackage[all]{xy}
\xyoption{arc} 

\theoremstyle{definition}

\makeatletter \@addtoreset{equation}{subsection} \makeatother  
\makeatletter \@addtoreset{equation}{subsubsection} \makeatother  

\DeclareFontFamily{OT1}{rsfs}{}
\DeclareFontShape{OT1}{rsfs}{n}{it}{<-> rsfs10}{}
\DeclareMathAlphabet{\curly}{OT1}{rsfs}{n}{it}

\newcommand{\Mod}{{\bf Mod}} 
\newcommand{\Vect}{{\bf Vect}} 
\newcommand{\Sets}{{\bf Ens}} 
\newcommand{\Ab}{{\bf A}}   
\newcommand{\B}{{\bf B}} 
\newcommand{\Sch}{{\bf Sch}} 

\newcommand{\EE}{\mathbb{E}} 
\newcommand{\GG}{\mathbb{G}} 
\newcommand{\LL}{\mathbb{L}} 

\newcommand{\slot}{ \hspace{0.05in} {\rm \_} \hspace{0.05in} } 

\newcommand{\bR}{\operatorname{\bf R}}
\newcommand{\bL}{\operatorname{\bf L}}

\renewcommand{\O}{\mathcal O} 

\newcommand{\ov}{\overline}
\newcommand\into{\hookrightarrow}

\newcommand{\Hom}{\operatorname{Hom}}   
\newcommand{\Ext}{\operatorname{Ext}}

\newcommand{\Tor}{\curly Tor}             
\newcommand{\sHom}{\curly Hom}            

\newcommand{\Ker}{\operatorname{Ker}}

\newcommand{\Sym}{\operatorname{Sym}}    

\newcommand{\Id}{\operatorname{Id}}
\newcommand{\Spec}{\operatorname{Spec}}

\newcommand{\Quot}{\operatorname{Quot}}

\newcommand{\at}{\operatorname{at}}

\newcommand{\Exal}{\operatorname{Exal}}

\begin{document}

\author{W.~D.~Gillam}
\address{Department of Mathematics, Brown University, Providence, RI 02912}
\email{wgillam@math.brown.edu}
\date{\today}
\title{Deformation of quotients on a product}

\begin{abstract} We consider the general problem of deforming a surjective map of modules $f : E \to F$ over a coproduct sheaf of rings $B=B_1 \otimes_A B_2$ when the domain module $E = B_1 \otimes_A E_2$ is obtained via extension of scalars from a $B_2$-module $E_2$.  Assuming $B_1$ is flat over $A$, we show that the Atiyah class morphism $F \to \LL_{B/B_2} \otimes^{\bL} F[1]$ in the derived category $D(B)$ factors naturally through (the shift of) a morphism $\beta : \Ker f \to \LL_{B/B_2} \otimes^{\bL} F$.  We describe the obstruction to lifting $f$ over a (square zero) extension $B_1' \to B_1$ in terms of $\beta$ and the class of the extension.  As an application, we use the reduced Atiyah class to construct a perfect obstruction theory on the Quot scheme of a vector bundle on a smooth curve (and more generally). \end{abstract}

\maketitle

\renewcommand{\theequation}{\arabic{equation}} 

\section*{Introduction} Let $Y$ be a scheme, $E$ a quasi-coherent sheaf on $Y$.  Let $X_0 \into X$ be a square zero closed immersion with ideal $I$ and let \begin{eqnarray} \label{SESa} 0 \to N_0 \to \pi_2^*E \to F_0 \to 0 \end{eqnarray} be an exact sequence of sheaves on $X_0 \times Y$ with $F_0$ flat over $X_0$.  Consider the problem of finding an exact sequence \begin{eqnarray} \label{SESb} 0 \to N \to \pi_2^*E \to F \to 0 \end{eqnarray} on $X \times Y$ with $F$ flat over $X$ restricting to \eqref{SESa} on $X_0 \times Y$.  Assuming $\Quot E$ is representable, this is equivalent to finding a map $X \to \Quot E$ restricting to the map $X_0 \to \Quot E$ determined by \eqref{SESa}.  It is ``well-known" that there is an obstruction $$\omega \in \Ext^1_{X_0 \times Y}(N_0, \pi_1^* I \otimes F_0) $$ whose vanishing is necessary and sufficient for the existence of such a sequence and that, when $\omega=0$, the set of such sequences is a torsor under $\Hom_{X_0 \times Y}(N_0, \pi_1^* I \otimes F_0)$.  A simple proof is provided for the reader's convenience in \eqref{elem}.

Following a trick of Nagata, Illusie (IV.3.2)\footnote{For a roman numeral $N$, a reference in the present paper of the form ($N.a$) always refers to section $a$ of chapter $N$ in Illusie's book \cite{Ill}.} constructs the obstruction $\omega$ by translating this lifting problem into a problem of $G := \GG_m$ equivariant algebra extensions, which can be handled using the machinery of the $G$ equivariant cotangent complex.  The basic link between the $G$ equivariant cotangent complex and the usual cotangent complex is provided by the equivalence of the two different descriptions of the Atiyah class in (IV.2.3.7.3).  This discussion motivates the main results of this paper, which we now briefly summarize.

For a scheme $X$ and a quasi-coherent sheaf $F$ on $X$, let $X[F]:=\Spec_X \O_X \oplus F$ denote the trivial square zero extension of $X$ by $F$.  The scheme $X[F]$ has a natural $G := \GG_m$ action obtained by restricting the action of scalars on $F$ to the invertible scalars.  $(X,F) \mapsto X[F]$ is a contravariant functor from the category (stack) of quasi-coherent sheaves to the category of $G$ schemes.  

Let $X,Y$ be schemes over a field $k$ so that the projection $\pi_2 : X \times Y \to Y$ is flat and we have a natural isomorphism $\LL_{X \times Y/Y} = \pi_1^* \LL_X$.  Let $E$ be a quasi-coherent sheaf on $Y$, \begin{eqnarray} \label{introSES} & \xymatrix{ 0 \ar[r] & N \ar[r] & \pi_2^*E \ar[r] & F \ar[r] & 0 } \end{eqnarray} an exact sequence of sheaves on $X \times Y$.  From the morphism of distinguished triangles of $G$ equivariant cotangent complexes associated to the diagram \begin{eqnarray} \label{Gschemes}  \xymatrix{ (X \times Y)[F] \ar[r] & (X \times Y)[\pi_2^*E] \ar[r] & Y[E] \\ (X \times Y)[F] \ar@{=}[u] \ar[r] & X \times Y \ar[r] \ar[u] & Y \ar[u] } \end{eqnarray} of $G$ schemes, we obtain a commutative diagram \begin{eqnarray} \label{redAtiyaha} \xymatrix{ \LL_{(X \times Y)[F]/(X \times Y)[\pi_2^*E]}^G \ar[r] & \LL_{(X \times Y)[\pi_2^*E]/Y[E]}^G \otimes^{\bL}_{\O_{X \times Y}[\pi_2^*E]} \O_{X \times Y}[F][1] \\ \LL_{(X \times Y)[F]/(X \times Y)}^G \ar[u] \ar[r] & \pi_1^*\LL_{X} \otimes_{\O_{X \times Y}}^{\bL} \O_{X \times Y}[F] \ar[u][1] }  \end{eqnarray} in the $G$ equivariant derived category $D^G((X \times Y)[F])$.  Pushing forward to $X \times Y$ and taking the weight one subcomplex is an exact functor \begin{eqnarray*} k^1 : D^G((X \times Y)[F]) \to D(X \times Y). \end{eqnarray*}  Applying $k^1$ to \eqref{redAtiyaha} we obtain a commutative diagram in $D(X \times Y)$ that can be written \begin{eqnarray} \label{redAtiyahb} \xymatrix@C+10pt{ N[1] \ar[r]^-{\beta[1]} & \pi_1^*\LL_{X} \otimes^{\bL} F[1] \\ F \ar[u] \ar[r]^-{{\rm at}(F)} & \pi_1^*\LL_{X} \otimes^{\bL} F \ar@{=}[u][1]. }  \end{eqnarray}  The bottom horizontal arrow is the Atiyah class of $F$ relative to $\pi_2$ and the $D(X \times Y)$ morphism $\beta : N \to \pi_1^*\LL_{X} \otimes^{\bL} F$ whose shift appears on the top line is called the \emph{reduced Atiyah class} of the quotient $\pi_2^*E \to F$.  It is studied extensively in Section~\ref{reducedatiyahclass} where we give a description of the reduced Atiyah class in terms of principal parts, thereby avoiding the equivariant cotangent complex.

Consider the functors \begin{eqnarray*} \Hom_{D(X)}( \LL_X, \slot), \; \Hom_{D(X \times Y)}(N, \pi_1^* \slot \otimes^{\bL} F) : D(X) \to \Vect_k \end{eqnarray*} and the ($k$ linear) natural transformation \begin{eqnarray*} \ov{\beta} : \Hom_{D(X)}( \LL_X, \slot) & \to & \Hom_{D(X \times Y)}(N, \pi_1^*  \slot \otimes^{\bL} F) \\ \ov{\beta}(M)(g) & := & (\pi_1^*g \otimes^{\bL} F) \beta. \end{eqnarray*}  Suppose that the Quot scheme of $E$ is representable, $X$ is an open subset of it, and \eqref{introSES} is the (restriction of the) universal sequence.  In Theorem~\ref{mainresult}, we show that, for any quasi-coherent sheaf $I$ on $X$ (regarded as a complex $I \in D(X)$ supported in degree zero), the $k$ vector space map \begin{eqnarray*} \ov{\beta}(I) : \Hom( \LL_X, I) & \to & \Hom(N, \pi_1^* I \otimes^{\bL} F) \end{eqnarray*} is an isomorphism and the map \begin{eqnarray*} \ov{\beta}(I[1]) : \Ext^1( \LL_X, I) & \to & \Ext^1(N, \pi_1^* I \otimes^{\bL} F) \end{eqnarray*} is injective.

If $Y$ is Gorenstein projective, then under some technical hypotheses, we show in Theorem~\ref{secondresult}, using Serre duality, that the functor $\Hom_{D(X \times Y)}(N, \pi_1^* F \otimes^{\bL} \slot)$ is represented by $$\EE := \bR \sHom( \bR \pi_{1*} \bR \sHom(N,F) , \O_X) \in D(X)$$ and that the $D(X)$ morphism $\EE \to \LL_X$ obtained from Yoneda's Lemma is an isomorphism on $H^0$ and surjective on $H^{-1}$.  If we assume furthermore that $N$ is locally free and $Y$ is a curve, then we will show that $\EE$ is of perfect amplitude $\subseteq [-1,0]$ and hence $\EE \to \LL_X$ is a perfect obstruction theory (POT) in the sense of \cite{BF}.  

Similar results have appeared in \cite{CFK} and \cite{MO}, but neither of these references explains the general mechanism for producing the derived category morphism $E \to \LL_Q$ ``intrinsically" using only the universal sequence on $X \times Y$.  In any case, the methods used here are completely different and should serve to clarify the construction of this POT and its variants, which are used in various places in the literature (\cite{Gil}, \cite{MO}, \cite{MOP}, \dots).

\subsection*{Acknowledgements}  Johan de Jong suggested the general idea of the reduced Atiyah class to me.  After thinking about it for a while, I realized that this was already implicit in (IV.3.2).  I owe a debt of gratitude to Luc Illusie since most of the results presented here are simple applications of the theory developed in his book \cite{Ill}.

\subsection*{Conventions.}  We work in a fixed topos $T$, as in Chapter IV, throughout; in all the applications $T$ is the topos of sheaves on a topological space.  We abbreviate ``ring object of $T$" to ``ring," etc.  When we speak of $H^i$ of a double (triple, etc.)\ complex, we mean $H^i$ of the associated total complex; we define ``quasi-isomorphism" of such complexes accordingly.  For a ring homomorphism $A \to B$, we write $I_\Delta \subseteq B \otimes_A B$ for the kernel of the multiplication map.  For a $B$ module $M$, let $$P^1_{B/A}(M)  :=  (B \otimes_A B / I_\Delta^2) \otimes_B M$$ where we use the right $B$ module structure (restriction of scalars along $b \mapsto 1 \otimes b$) on $(B \otimes_A B / I_\Delta^2)$ to define the tensor product over $B$, then we regard the result as a $B$ module via the left $B$ module structure (restriction of scalars along $b \mapsto b \otimes 1$).  ``Extension" always means ``square zero extension" (surjective map of algebras whose kernel is a square zero ideal).  We use $k^n$ to denote the functor from graded modules (resp.\ complexes of graded modules, the derived category of the category of graded modules, \dots) over a graded ring $B$ to modules over the degree zero part of $B$ given by taking the degree $n$ part of a module (resp.\ complex  of modules; this is exact functor).  A \emph{retract} of a morphism $f : A \to B$ (in any category) is a morphism $s : B \to A$ such that $sf = \Id_A$, while a \emph{section} of $f$ is a morphism $ s : B \to A$ such that $fs = \Id_B$.  These are dual: $s$ is a retract of $f$ iff $s^{\rm op} : A \to B$ is a section of $f^{\rm op} : B \to A$ in the opposite category.

\section{Deformation of quotients on a coproduct}

\renewcommand{\theequation}{\thesubsection.\arabic{equation}} 

\subsection{Setup} \label{setup} The following basic setup will be used throughout this paper.  Let $A \to B_0$ be a ring homomorphism and assume $B_0=B_1 \otimes_A B_2$ is the coproduct of $A$ algebras $B_1,B_2$.  Let $E_2$ be a $B_2$ module and let $$E_0 := B_0 \otimes_{B_2} E_2 = B_1 \otimes_A E_2$$ be the $B_0$ module obtained from $E_2$ by extension of scalars.  Let \begin{eqnarray} \label{SES2} \xymatrix{ 0  \ar[r] & N_0 \ar[r] & E_0 \ar[r]^{f_0} & F_0 \ar[r] & 0 } \end{eqnarray} be an exact sequence of $B_0$ modules and let \begin{eqnarray} \label{extension1} \underline{B}_1 & = & \xymatrix{ 0 \ar[r] & I_1 \ar[r] & B_1' \ar[r] & B_1 \ar[r] & 0 \\ & & & A \ar[ul] \ar[u] } \end{eqnarray} be an $A$ extension of $B_1$ by a $B_1$ module $I_1$.  Assume that \begin{eqnarray} \label{tor1vanishes1} \Tor_1^A(B_1,B_2) &=& 0. \end{eqnarray} Let $I=I_1 \otimes_A B_2$ and let \begin{eqnarray} \label{extension2} \underline{B} & = & \xymatrix{ 0 \ar[r] & I \ar[r] & B \ar[r] & B_0 \ar[r] & 0 \\ & & & A \ar[ul] \ar[u] } \end{eqnarray} be the $A$ extension of $B_0$ by $I$ obtained from \eqref{extension1} by applying $\otimes_A B_2$.

Let $E = B_1' \otimes_A B_2$ and assume that  \begin{eqnarray} \label{tor1vanishes2} \Tor_1^A(B_1,E_2) &=& 0 \end{eqnarray} so we have an exact sequence \begin{eqnarray} \label{SES3} \underline{E} = ( \xymatrix{ 0  \ar[r] & I \otimes_{B_0} E_0 \ar[r] & E \ar[r] & E_0 \ar[r] & 0 }) \end{eqnarray} of $B$ modules obtained from \eqref{extension1} by taking $\otimes_A E_2$.  (As usual, we identify $B_0$ modules with $B$ modules annihilated by $I$.)

\subsection{Problem} \label{liftingproblem} Fix a $B_0$ module $K$ and a morphism $u_0 : I \otimes_{B_0} F_0 \to K$ of $B_0$ modules. Consider the problem of finding a $B$ module extension $\underline{F}$ of $F_0$ by $K$ and a morphism of extensions $\underline{f} : \underline{E} \to \underline{F}$ of the form: \begin{eqnarray*} \label{liftingproblem2} \xymatrix{ 0  \ar[r] & I \otimes_{B_0} E_0 \ar[r] \ar[d]_{u_0(I \otimes f_0)} & E \ar[d] \ar[r] & E_0 \ar[r] \ar[d]^{f_0} & 0 \\ 0 \ar[r] & K \ar[r] & F \ar[r] & F_0 \ar[r] & 0 } \end{eqnarray*}

The basic result concerning this problem is the following ``classical" theorem, which we will prove in \eqref{elem} and again in \eqref{Illusie}.

\subsection{Theorem} \label{mainthm} \emph{There is an obstruction $\omega \in \Ext^1_{B_0}(N_0, K)$ whose vanishing is necessary and sufficient for the existence of a solution to \eqref{liftingproblem}.  When $\omega=0$, the set of solutions to \eqref{liftingproblem} is a torsor under $\Hom_{B_0}(N_0,K)$.}

\subsection{Remark} An extreme case of the problem \eqref{liftingproblem} occurs when $u_0=0$ (equivalently $u_0(I \otimes f_0)=0$).  In this case $f_0$ always factors through $F \to F_0$ and we see that a solution to the problem is the same thing as a morphism of extensions \begin{eqnarray} \label{liftingproblem3} \xymatrix{ 0 \ar[r] & N_0 \ar[r] \ar[d] & E_0 \ar[d] \ar[r] & F_0 \ar[r] \ar@{=}[d] & 0 \\ 0 \ar[r] & K \ar[r] & F \ar[r] & F_0 \ar[r] & 0. } \end{eqnarray} Of course, every such morphism is obtained by pushing out the top row along a unique morphism $N_0 \to K$ so in this case it is clear that the obstruction vanishes and, in fact, we have a natural trivialization of the torsor under $\Hom(N_0,K)$ because we have a natural base point obtained from the split extension.

\subsection{Remark} \label{flatdeformations} Assume $F_0$ is flat over $B_1$.  Then by the flatness criterion, $F$ will be flat over $B_1'$ iff the natural map $u(\underline{F}) : I \otimes_{B_0} F_0 \to K$ (IV.3.1) is an isomorphism (note that $I \otimes_{B_0} F_0 = I_1 \otimes_{B_1} F_0$).  Using this isomorphism to make the identification $I \otimes_{B_0} F_0 = K$, we see that the problem of lifting $f_0 : E_0 \to F_0$ to a surjection $f : E \to F$ with $F$ flat over $B_1'$ is equivalent to the special case of \eqref{liftingproblem} where $K = I \otimes_{B_0} F_0$ and $u_0$ is the identity.  In this case, we refer to a solution of \eqref{liftingproblem} as a \emph{flat deformation} of $f_0$ over $B$.

\subsection{} \label{elem} The obstruction $\omega \in \Ext^1_{B_0}(N_0, K)$ can be constructed elementarily as follows.  Let $L$ be the kernel of the composition $E \to E_0 \to F_0$.  We have a commutative diagram \begin{eqnarray} \label{elem1} \xymatrix{ & 0 \ar[d] & 0 \ar[d] & 0 \ar[d] \\ 0 \ar[r] & I \otimes_{B_0} E_0 \ar[r] \ar@{=}[d] & L \ar[d] \ar[r] & N_0 \ar[d] \ar[r] & 0 \\ 0 \ar[r] & I \otimes_{B_0} E_0 \ar[r] \ar[d] & E \ar[r] \ar[d] & E_0 \ar[r] \ar[d]^{f_0} & 0 \\ & 0 \ar[r] & F_0 \ar[d] \ar@{=}[r] & F_0 \ar[d] \ar[r] & 0 \\ & & 0 & 0 } \end{eqnarray} with exact columns and rows.  By pushing out the top row along $u_0(I \otimes f_0)$ we obtain a morphism of extensions \begin{eqnarray} \label{elem2} \xymatrix{ 0 \ar[r] & I \otimes_{B_0} E_0 \ar[d]_{u_0(I \otimes f_0)} \ar[r] & L \ar[r] \ar[d] & N_0 \ar[r] \ar@{=}[d] & 0 \\ 0 \ar[r] & K \ar[r] & M \ar[r]  & N_0 \ar[r] & 0 } \end{eqnarray} where $M:=K \oplus_{I \otimes E_0} L$.  Note that $M$ is a $B$ module (i.e.\ $IM=0$): $$[ik,il]=[0,il]=[i f_0(l),0]=0,$$ so the bottom row defines a $B_0$ module extension.  Let $\omega \in \Ext^1_{B_0}(N_0, K)$ be the class of this extension.   I claim that the vanishing of $\omega$ is necessary and sufficient for the existence of a solution to \eqref{liftingproblem}, and that, in fact, splittings of the bottom row in \eqref{elem2} (which are a pseudo-torsor under $\Hom_{B_0}(N_0,K)$) correspond bijectively to solutions of \eqref{liftingproblem}.

Given a solution $\underline{f} : \underline{E} \to \underline{F}$ to \eqref{liftingproblem} we obtain a commutative diagram \begin{eqnarray} \label{elem5} \xymatrix{ 0 \ar[r] & I \otimes_{B_0} E_0 \ar[r] \ar@{=}[d] & N \ar[d] \ar[r] & N_0 \ar[r] \ar[d] & 0 \\ 0 \ar[r] & I \otimes_{B_0} E_0 \ar[r] \ar[d]_{u_0(I \otimes f_0)} \ar[r] & E \ar[r] \ar[d]^f & E_0 \ar[d]^{f_0} \ar[r] & 0 \\ 0 \ar[r] & K \ar[r] & F \ar[r] & F_0 \ar[r] & 0 }  \end{eqnarray} with exact rows by pulling back the extension in the middle row along $N_0 \to E_0$.  Note that $N \to E$ is monic and $N \subset L \subset E$, but $N$ need not be contained in the kernel of $f$.  Given $n \in N_0$, choose a (local) lift $\ov{n} \in N$.  I claim that $n \mapsto [-f(\ov{n}),\ov{n}] \in M$ is well-defined (independent of the choice of lift $\ov{n} \in N$), in which case it will clearly provide a splitting of the bottom row in \eqref{elem2}.  Indeed, if $\ov{n}'$ is another (local) lift, then $\ov{n}'=\ov{n} + \epsilon$ for some $\epsilon \in I \otimes_{B_0} E_0$ and we have $$[-f(\ov{n} + \epsilon),\ov{n} + \epsilon] = [-f(\ov{n})-u_0(I \otimes f_0)(\epsilon), \ov{n}+\epsilon] = [-f(\ov{n}),\ov{n}] \in M.$$

Conversely, notice that a splitting of the bottom row in \eqref{elem2} is the same thing as a map $s : L \to K$ making the diagram \begin{eqnarray} \label{elem3} \xymatrix{ I \otimes_{B_0} E_0 \ar[d]_{u_0(I \otimes f_0)} \ar[r] & L \ar[ld]^s \\ K } \end{eqnarray} commute.  By pushing out the sequence defining $L$ along $s$ we obtain a commutative diagram \begin{eqnarray} \label{elem4} \xymatrix{ 0 \ar[r] & I \otimes_{B_0} E_0 \ar[r] \ar[d] & E \ar@{=}[d] \ar[r] & E_0 \ar[d]^{f_0} \ar[r] & 0 \\ 0 \ar[r] & L \ar[d]^-s \ar[r] & E \ar[d] \ar[r] & F_0 \ar@{=}[d] \ar[r] & 0 \\ 0 \ar[r] & K \ar[r] & F \ar[r] & F_0 \ar[r] & 0 } \end{eqnarray} with exact rows, where we have set $F := K \oplus_L E$.  The map from the top row to the bottom row is evidently a solution to \eqref{liftingproblem}.

It is straightforward to check that the two constructions are inverse.

\subsection{Remark} \label{torsorbijection}  Assume $u_0$ is surjective, hence $u_0(I \otimes f_0)$ is also surjective.  Let $S$ be its kernel.  Given any solution $\underline{f} : \underline{E} \to \underline{F}$ to \eqref{liftingproblem}, set $N := \Ker f$.  We have a commutative diagram \begin{eqnarray} & \xymatrix{ & 0 \ar[d] & 0 \ar[d] & 0 \ar[d] \\ 0 \ar[r] & S \ar[d] \ar[r] & N \ar[d] \ar[r] & N_0 \ar[d] \ar[r] & 0   \\ 0 \ar[r] & I \otimes_{B_0} E_0 \ar[d] \ar[r] & E \ar[d] \ar[r] & E_0 \ar[d] \ar[r] & 0 \\ 0 \ar[r] & K \ar[r] \ar[d] & F \ar[r] \ar[d] & F_0 \ar[r] \ar[d] & 0 \\ & 0 & 0 & 0 } \label{notation} \end{eqnarray} with exact rows and columns.  Suppose $\underline{f}^* : \underline{E} \to \underline{F}^*$ is some fixed solution to \eqref{liftingproblem}.  Then by following through the above discussion we see that the bijection between solutions to \eqref{liftingproblem} and $\Hom_{B_0}(N_0,K)$ (using $\underline{f}^*$ as a basepoint to trivialize the torsor) is realized by \begin{eqnarray} ( \underline{f} : \underline{E} \to \underline{F} ) & \mapsto &  \xymatrix{ S \ar[r] & N \ar[d] \\ & F^* \ar[r] & F_0. }  \end{eqnarray}  Here, the top row on the right is a complex quasi-isomorphic to $N_0$ and the bottom row on the right is a complex quasi-isomorphic to $K$.  The map $N \to F^*$ is the composition of the inclusion $N \to E$ and $f : E \to F^*$.

\subsection{} \label{Illusie} Illusie constructs the obstruction $\omega$, following a trick of Nagata, by translating the problem \eqref{liftingproblem} into a problem of graded algebra extensions and then using the machinery of the graded cotangent complex.  Let $$B[E] := \Sym^*_{B} E / \Sym^{>1}_B E$$ be the algebra of dual numbers on $E$ over $B$, graded as usual so that $E$ is the degree one component; define graded algebras $B_0[E_0]$ and $B_0[F_0]$ similarly.  There is an obvious diagram \begin{eqnarray} \label{gradedalgebras} \xymatrix{ B[E] \ar[r] & B_0[E_0] \ar[r] & B_0[F_0] } \end{eqnarray} of graded rings.  Regard $K$ as a graded $B_0[F_0]$ module supported in degree one (hence annihilated by $F_0$).  A solution to \eqref{liftingproblem} is the same thing as an extension of graded algebras \begin{eqnarray} \label{extensiongr} \xymatrix{ 0 \ar[r] & K \ar[r] & G \ar[r] & B_0[F_0] \ar[r] & 0 \\ & & & B[E] \ar[ul] \ar[u] } \end{eqnarray} where induced map $I \otimes_{B_0} E_0 \to E \to K$ is $u_0(I \otimes f_0)$.  Indeed, the bijection is given by taking the degree one part of $G$ (the extension \eqref{extensiongr} is uninteresting in all other degrees since $K$ is supported in degree one, so we can always write $G=B_0[F]$ for some $B$ module $F$).

Graded algebra extensions \eqref{extensiongr} form a graded $B_0[E_0]$ module $\Exal_{B[E]}^{\rm gr}(B_0[F_0],K)$.  By the fundamental theorem of the graded cotangent complex (IV.2.4.2) and (IV.1.2.2.1) there are isomorphisms \begin{eqnarray} \label{exaliso} \Exal_{B[E]}^{\rm gr}(B_0[F_0],K) & = & \Ext^1_{B_0[F_0]}(\LL_{B_0[F_0]/B[E]}^{\rm gr} , K)_{\rm gr} \\ \nonumber &=& \Ext^1_{B_0}(k^1 \LL_{B_0[F_0]/B[E]}^{\rm gr} , K). \end{eqnarray}  From (IV.2.2.5) and the triangle associated to $B_0 \to B_0[E_0] \to B_0[F_0]$ we obtain a natural isomorphism \begin{eqnarray} \label{iso} k^1 \LL_{B_0[F_0]/B_0[E_0]}^{\rm gr} = N_0[1]. \end{eqnarray}  Using \eqref{iso}, the degree one part of the distinguished triangle of cotangent complexes associated to \eqref{gradedalgebras} is a triangle in $D(B_0)$ that can be written \begin{eqnarray} \label{k1triangle} \xymatrix{ k^1 \LL_{B_0[E_0]/B[E]}^{\rm gr} \otimes_{B_0[E_0]}^{\bL} B_0[F_0] \ar[r] & k^1 \LL_{B_0[F_0]/B[E]}^{\rm gr}  \ar[r] & N_0[1]  }. \end{eqnarray}  Since $B[E] \to B_0[E_0]$ is a square zero extension with kernel $I \oplus (I \otimes_{B_0} E_0)$ (in gradings $0,1$) we have $$\tau_{\geq -1} \LL_{B_0[E_0]/B[E]}^{\rm gr} = (I \oplus (I \otimes_{B_0} E_0))[1]$$ and hence \begin{eqnarray} \label{iso2} \tau_{\geq -1} \LL_{B_0[E_0]/B[E]}^{\rm gr} \otimes^{\bL}_{B_0[E_0]} B_0[F_0] &=& (I \oplus (I \otimes_{B_0} F_0))[1] \end{eqnarray} (in gradings $0,1$).   Applying $\Hom( \slot , K)$ to \eqref{k1triangle} and using \eqref{exaliso} and \eqref{iso2} we obtain a long exact sequence \begin{eqnarray} \label{LES} & \xymatrix@R-15pt{ & 0 \ar[r] & \Hom(N_0,K) \ar[r] & \Exal_{B[E]}^{\rm gr}(B_0[F_0],K) \ar[r] & \Hom(I \otimes_{B_0} F_0,K) \\ & \ar[r]^-{\delta} & \Ext^1(N_0,K) & \cdots } \end{eqnarray} of $B_0$ modules.  A solution to \eqref{liftingproblem} is an element of $\Exal_{B[E]}^{\rm gr}(B_0[F_0],K)$ whose image in $\Hom(I \otimes_{B_0} F_0,K)$ is $u_0$.  Evidently such a solution exists iff $\omega := \delta u_0$ vanishes, in which case solutions are a torsor under $\Hom(N_0,K)$.

As a morphism in $D(B_0)$, $\omega$ is the composition \begin{eqnarray} \label{omega} & \xymatrix{ N_0 \ar[r] & k^1 \LL_{B_0[E_0]/B[E]} \otimes^{\bL}_{B_0[E_0]} B_0[F_0] \ar[r] & I \otimes_{B_0} F_0[1] \ar[r]^-{u_0[1]} & K[1] } \end{eqnarray} where the first map is the one from the triangle \eqref{k1triangle} and the second is the truncation.  

\subsection{} As the notation suggests, the obstruction $\omega$ constructed in \eqref{elem} coincides with the one constructed in \eqref{Illusie}.  To see this, first recall that we natural isomorphisms $$k^1 \LL_{B_0[F_0]/B_0[E_0]} = \tau_{\geq -1} k^1 \LL_{B_0[F_0]/B_0[E_0]} = N_0[1],$$ so the morphism $N_0 \to I \otimes_{B_0} F_0[1]$ constructed in \eqref{Illusie} using the degree one part of the transitivity triangle associated to \eqref{gradedalgebras} coincides with the morphism $$ \tau_{\geq -1} k^1 \LL_{B_0[F_0]/B_0[N_0]}[-1] \to \tau_{\geq -1} k^1 \LL_{B_0[E_0]/B[E]} \otimes^{\bL}_{B_0[E_0]} B_0[F_0] = I \otimes_{B_0} F_0[1] $$ obtained from the degree one part of the truncated transitivity triangle \begin{eqnarray} \label{truncatedtriangle} & \xymatrix@C-8pt{ \tau_{\geq -1} \LL_{B_0[E_0]/B[E]} \otimes^{\bL}_{B_0[E_0]} B_0[F_0] \ar[r] & \tau_{\geq -1} \LL_{B_0[F_0]/B[E]} \ar[r] & \tau_{\geq -1} \LL_{B_0[F_0]/B_0[E_0]}. } \end{eqnarray} 

The graded algebra maps appearing in \eqref{gradedalgebras} are all surjective, with kernels as indicated below.  $$ \xymatrix{ B_0[E_0] \ar[r]^-{0 \oplus N_0} & B_0[F_0] \\ B[E] \ar[u]^{I \oplus (I \otimes_{B_0} E_0)} \ar[ru]_-{I \oplus L} }$$  In each case, the direct summand decomposition corresponds to the splitting into the degree $0,1$ parts, respectively.  Recall that \begin{eqnarray*} ( I \oplus (I \otimes_{B_0} E_0) ) \otimes_{B_0[E_0]} B_0[F_0] &=& \displaystyle{ \frac{I \oplus (I \otimes_{B_0} E_0)}{N_0 \cdot (I \oplus (I \otimes_{B_0} E_0)) } } \\ &=& \displaystyle{ \frac{I \oplus (I \otimes_{B_0} E_0)}{0 \oplus (I \otimes_{B_0} N_0) } } \\ &=& I \oplus (I \otimes_{B_0} F_0). \end{eqnarray*} By standard facts about the truncated cotangent complex (III.1.3), the triangle \eqref{truncatedtriangle} is naturally identified, after applying $[-1]$, with (the triangle associated to) the short exact sequence \begin{eqnarray} \label{seq} & \xymatrix{ 0 \ar[r] & I  \oplus (I \otimes_{B_0} F_0) \ar[r] & \displaystyle{\frac{I \oplus L}{(I \oplus L)^2}} \ar[r] & 0 \oplus N_0 \ar[r] & 0 } \end{eqnarray} of graded $B_0[F_0]$ modules.  Note that $I \oplus (I \otimes_{B_0} E_0) \subset B[E]$ and $0 \oplus N_0 \subseteq B_0[E_0]$ are square zero ideals, while $$(I \oplus L)^2 = 0 \oplus IL \subseteq I \oplus L.$$  The degree one part of the exact sequence \eqref{seq} is therefore \begin{eqnarray} \label{diag2} & \xymatrix{ 0 \ar[r] & I \otimes_{B_0} F_0 \ar[r] & L/IL \ar[r] &  N_0 \ar[r] & 0. } \end{eqnarray}

There is a natural isomorphism $I \otimes_{B_0} N_0 = IL$ obtained by mapping $i \otimes n$ to $il$ where $l$ is a (local) lift of $n$ under the surjection $L \to N_0$ (the choice of local lift is irrelavant since the kernel $I \otimes_{B_0} E_0$ of $L \to N_0$ is annihilated by $I$).  In fact, \eqref{diag2} is part of the commuatative diagram \begin{eqnarray} \label{diag} & \xymatrix{ & 0 \ar[d] & 0 \ar[d] \\ 0 \ar[r] & I \otimes_{B_0} N_0 \ar@{=}[r] \ar[d] & IL \ar[d] \ar[r] & 0 \ar[d] \\ 0 \ar[r] & I \otimes_{B_0} E_0 \ar[d]_{I \otimes f_0} \ar[r] & L \ar[d] \ar[r] & N_0 \ar@{=}[d] \ar[r] & 0 \\ 0 \ar[r] & I \otimes_{B_0} F_0 \ar[d] \ar[r] & L/IL \ar[d] \ar[r] & N_0 \ar[r] \ar[d] & 0 \\ & 0 & 0 & 0} \end{eqnarray} with exact columns and rows.  

As in any map of extensions with isomorphic cokernels, the bottom left square in \eqref{diag} is a pushout.  Therefore, upon pushing out the bottom row along $u_0 : I \otimes_{B_0} F_0 \to K$, the resulting extension of $N_0$ by $K$ coincides with the one on the bottom row of \eqref{elem2}.  Note that $L/IL$ is manifestly a $B_0$ module (annihilated by $I$), so this gives another proof that the bottom row of \eqref{elem2} is a $B_0$ module extension.  On the other hand, the map $N_0 \to K[1]$ obtained from this extension coincides with the composition of the map $N_0 \to I \otimes_{B_0} F_0[1]$ (obtained from the truncated triangle \eqref{truncatedtriangle}) and $u_0[1]$.

\subsection{} So far we have only reviewed a special case of (IV.3.2).  We have not yet made use of the fact that $B_0 = B_1 \otimes_A B_2$ is a coproduct or the fact that the extension $\underline{B}$ \eqref{extension2} of $B_0$ is obtained by base change from the extension $\underline{B}_1$ \eqref{extension1} of $B_1$.  The rest of this section is devoted to giving an alternative description of the morphism $\omega$, assuming $B_1$ is flat over $A$, in terms of the class \begin{eqnarray} \label{eB1} e(\underline{B}_1) \in \Ext^1_{B_1}(\LL_{B_1/A},I_1) \end{eqnarray} corresponding to the extension \eqref{extension1} under the isomorphism provided by the fundamental theorem of the cotangent complex (III.1.2.3) and the ``reduced Atiyah class" of the quotient map $f_0$, to be introduced below.

\subsection{Assumption} \label{assumptions} We assume throughout that $B_1$ is flat over $A$, and hence $B_0$ is flat over $B_2$.  Note that this assumption implies the assumptions \eqref{tor1vanishes1} and \eqref{tor1vanishes2} of \eqref{setup}.  It also implies that $B_1$ and $B_2$ are \emph{tor-independent} over $A$, in the sense that $\Tor^A_{>0}(B_1,B_2)=0$.  It is possible to get away with significantly less, but the concomitant complications would make the exposition even more obtuse than it already is.  In the applications of Section~\ref{applications}, the ring $A$ will have tor dimension zero, hence $B_1$ will be trivially $A$ flat.

\subsection{} \label{atiyahclassintro} By taking the degree one part $k^1$ of the map of cotangent complex transitivity triangles associated to \begin{eqnarray} \label{source1} & \xymatrix{ B[E] \ar[r] & B_0[E_0] \ar[r] & B_0[F_0] \\ B_2[E_2] \ar[u] \ar[r] & B_0[E_0] \ar@{=}[u] \ar[r] & B_0[F_0] \ar@{=}[u] } \end{eqnarray} we obtain a commutative diagram \begin{eqnarray} \label{diagram1} \xymatrix@C+15pt{ N_0 \ar[r] \ar@{=}[d]  & k^1 \LL^{\rm gr}_{B_0[E_0]/B[E]} \otimes_{B_0[E_0]}^{\bL} B_0[F_0] \\ N_0 \ar[r] & k^1 \LL_{B_0[E_0]/B_2[E_2]}^{\rm gr} \otimes_{B_0[E_0]}^{\bL} B_0[F_0] \ar[u] }  \end{eqnarray} in $D(B_0)$.  

From the commutative diagram of graded rings \begin{eqnarray} \label{source2} \xymatrix{ & B \ar[rd] \ar[dd] \\ B_2 \ar[ru] \ar[rr] \ar[dd] & & B_0 \ar[dd] \\ & B[E] \ar[rd] \\ B_2[E_2] \ar[rr] \ar[ru] & & B_0[E_0] } \end{eqnarray} and the naturality of the (graded) cotangent complex, we obtain a commutative diagram \begin{eqnarray} \label{diagram2} \xymatrix{ \LL_{B_0/B} \otimes_{B_0}^{\bL} B_0[E_0] \ar[r] & \LL_{B_0[E_0]/B[E]}^{\rm gr} \\ \LL_{B_0/B_2} \otimes_{B_0}^{\bL} B_0[E_0] \ar[u] \ar[r] & \LL_{B_0[E_0]/B_2[E_2]}^{\rm gr} \ar[u] } \end{eqnarray} in $D(B_0[E_0])_{\rm gr}$.  Since the front square in \eqref{source2} is cocartesian: $$B_0[E_0]=B_0 \otimes_{B_2} B_2[E_2],$$ it follows from the base change properties of the cotangent complex (II.2.2) and the fact that $B_0$ is flat over $B_2$ \eqref{assumptions} that the bottom arrow in \eqref{diagram2} is an isomorphism.\footnote{All that is really needed is that $\Tor^{B_2}_{>0}(B_0,E_2)=0$.}  By applying $\otimes_{B_0[E_0]}^{\bL} B_0[F_0]$ to \eqref{diagram2} and taking $k^1$, we obtain a commutative diagram \begin{eqnarray} \label{diagram3} \xymatrix{ \LL_{B_0/B} \otimes_{B_0}^{\bL} F_0 \ar[r] & k^1 \LL_{B_0[E_0]/B[E]}^{\rm gr} \otimes_{B_0[E_0]}^{\bL} B_0[F_0] \\ \LL_{B_0/B_2} \otimes_{B_0}^{\bL} F_0 \ar[u] \ar[r] & k^1 \LL_{B_0[E_0]/B_2[E_2]}^{\rm gr} \otimes_{B_0[E_0]}^{\bL} B_0[F_0] \ar[u] } \end{eqnarray} in $D(B_0)$, where the bottom arrow is an isomorphism.

From the maps of cotangent complex transitivity triangles associated to the diagram of rings \begin{eqnarray} \label{source3} \xymatrix{ A \ar@{=}[d] \ar[r] & B_1' \ar[d] \ar[r] & B_1 \ar[d] \\ A \ar[r] & B \ar[r] & B_0 \\ A \ar@{=}[u] \ar[r] & B_2 \ar[r] \ar[u] & B_0 \ar@{=}[u] } \end{eqnarray} we obtain a commutative diagram \begin{eqnarray} \label{diagram4} \xymatrix{ \LL_{B_1/A} \otimes^{\bL}_{B_1} B_0 \ar[r] \ar[d] & \LL_{B_1/B_1'} \otimes^{\bL}_{B_1} B_0 \ar[r] \ar[d] & I_1 \otimes_{B_1} B_0[1] \ar@{=}[d] \\ \LL_{B_0/A} \ar[r] \ar[rd] & \LL_{B_0/B} \ar[r] & I[1] \\ & \LL_{B_0/B_2} \ar[u] } \end{eqnarray} in $D(B_0)$.  The two right horizontal arrows are the truncations $\tau_{\geq -1}$ and the composition given by the top row is $e(\underline{B}_1) \otimes_{B_1} B_0$, where $$e(\underline{B}_1) \in \Ext^1(\LL_{B_1/A},I_1)$$ corresponds to the algebra extension $\underline{B}_1$ \eqref{extension1} under the Fundamental Theorem (III.1.2.3) .  Under the tor-independence assumption \eqref{assumptions}, the base change properties of the cotangent complex (II.2.2) imply that the composition $\LL_{B_1/A} \otimes_{B_1}^{\bL} B_0 \to \LL_{B_0/B_2}$ in this diagram is an isomorphism.  Combining this isomorphism with the one given by the bottom row of \eqref{diagram3} yields a natural isomorphism \begin{eqnarray} k^1 \LL_{B_0[E_0]/B_2[E_2]}^{\rm gr} \otimes_{B_0[E_0]}^{\bL} B_0[F_0] & = & \LL_{B_1/A} \otimes_{B_1}^{\bL} F_0, \end{eqnarray} hence the bottom arrow in \eqref{diagram1} may be viewed as a morphism \begin{eqnarray} \label{reducedatiyah} {\rm ra}(f_0) : N_0 & \to & \LL_{B_1/A} \otimes_{B_1}^{\bL} F_0 \end{eqnarray} which we will call the \emph{reduced Atiyah class} of the quotient $f_0$.  It will be further studied in Section~\ref{reducedatiyahclass}.

By assembling the diagrams \eqref{diagram1}, \eqref{diagram3}, and \eqref{diagram4}$\otimes_{B_0}^{\bL} F_0$ and using the isomorphisms noted above, we obtain a commutative diagram \begin{eqnarray} \label{diagram5} \xymatrix{ N_0 \ar[r] \ar[rd]_-{{\rm ra}(f_0)} & k^1 \LL_{B_0[E_0]/B[E]}^{\rm gr} \otimes^{\bL}_{B_0[E_0]} B_0[F_0] \ar[r] & I \otimes_{B_0}^{\bL} F_0[1] \\ & \LL_{B_1/A} \otimes_{B_1}^{\bL} F_0 \ar[ru]_-{e(\underline{B}_1) \otimes_{B_1} F_0} } \end{eqnarray} in $D(B_0)$.  The composition of the top row and $u_0[1]$ is the obstruction $\omega$ \eqref{omega}. We have proved the following (c.f.\ (IV.3.1.8), (IV.3.2.14), and Theorem~2.9 in \cite{HT}):

\subsection{Theorem}  \label{atiyahclassthm} \emph{ Assume $B_1$ is flat over $A$.  Then the obstruction $\omega \in \Ext^1(N_0,K)$ of \eqref{mainthm} can be written as the composition \begin{eqnarray*} & \xymatrix@C+10pt{ N_0 \ar[r] & \LL_{B_1/A} \otimes_{B_1}^{\bL} F_0 \ar[r] & I \otimes_{B_0} F_0[1] \ar[r]^-{u_0[1]} & K[1], } \end{eqnarray*} where the first arrow is given by the reduced Atiyah class of $f_0$ \eqref{reducedatiyah} and the second is $e(\underline{B}_1) \otimes_{B_1} F_0$, where $e(\underline{B}_1) \in \Ext^1(\LL_{B_1/A},I_1)$ corresponds to $\underline{B}_1$ under the Fundamental Theorem (III.1.2.3).}

\subsection{} For later use, we will now determine the map induced on $H^0$ by the reduced Atiyah class $N_0 \to \LL_{B_1/A} \otimes^{\bL}_{B_1} F_0$ \eqref{reducedatiyah}.  Recall that this map was constructed from the natural diagram of solid arrows \begin{eqnarray} & \label{di} \xymatrix{ N_0  \ar@{.>}[rdd]_-{{\rm ra}(f_0)} \ar[r] & k^1 \LL_{B_0[E_0]/B_2[E_2]}^{\rm gr} \otimes^{\bL}_{ B_0[E_0]} B_0[F_0] \\ & \LL_{B_0/B_2} \otimes_{B_0}^{\bL} F_0 \ar[u]_{\cong} \\ & \LL_{B_1/A} \otimes^{\bL}_{B_1} F_0 \ar[u]_{\cong} } \end{eqnarray} and the fact that the indicated arrows are isomorphisms under the assumption \eqref{assumptions} that $B_1$ is flat over $A$.  These isomorphisms were obtained from the base change property (II.2.2) of the cotangent complex.  Note that, even without the flatness assumption \eqref{assumptions} the indicated arrows induce isomorphisms on $H^0$ because K\"ahler differentials commute with all direct limits, in particular with pushouts.  

Taking $H^0$ of \eqref{di}, we obtain the diagram: \begin{eqnarray} & \label{H0di} \xymatrix{ N_0 \ar[r] & k^1\Omega_{B_0[E_0]/B_2[E_2]} \otimes_{ B_0[E_0]} B_0[F_0] \\ & \Omega_{B_0/B_2} \otimes_{B_0} F_0 \ar@{=}[u] \\ & \Omega_{B_1/A} \otimes_{B_1} F_0 \ar@{=}[u] } \end{eqnarray} The horizontal arrow is just the ``connecting homomorphism" from the long exact sequence obtained from the triangle associated to the bottom row of \eqref{source1} (note $N_0 = H^{-1}(\LL_{B_0[F_0]/B_0[E_0]}^{\rm gr})$).  Viewing $N_0$ as a submodule of $E_0 = B_1 \otimes_A E_2$, this map is given by $b_1 \otimes e \mapsto d(b_1 \otimes e) \otimes 1$.  Note $b_1 \otimes e$ is in the degree one part of $B_0[E_0]$, so $d(b_1 \otimes e)$ is in $k^1 \Omega_{B_0[E_0]/B_0[F_0]}$.  Notice that \begin{eqnarray*} d(b_1 \otimes e) \otimes 1 &=& d((1 \otimes e)(b_1 \otimes 1)) \otimes 1 \\ &=& (1 \otimes e) \cdot d(b_1 \otimes 1) \otimes 1 \\ &=& d(b_1 \otimes 1) \otimes f_0(1 \otimes e) \end{eqnarray*} is the image of $db_1 \otimes f_0(1 \otimes e)$ under the vertical isomorphisms, so the map induced on $H^0$ by the reduced Atiyah class is given by \begin{eqnarray} \label{H0ra} H^0({\rm ra}(f_0)) : N_0 & \to & \Omega_{B_1/A} \otimes_{B_1} F_0 \\ \nonumber b_1 \otimes e & \mapsto & db_1 \otimes f_0(1 \otimes e). \end{eqnarray}  Notice that this map is $B_0$ linear because $b_1 \otimes e \in N_0 \subset E_0 = B_1 \otimes_A E_2$.

\section{Reduced Atiyah class}  \label{reducedatiyahclass} In this section, we give an independent treatment of the reduced Atiyah class, first introduced in \eqref{reducedatiyah} above.  The results of this section are not needed elsewhere; our purpose here is only to further examine the reduced Atiyah class and give an alternative construction of it via principal parts.

\subsection{} \label{redatiyah1} Let $B \to C$ be a ring homomorphism and let $E$ be a $C$ module.  Recall (III.1.2.6.3) the principal parts sequence \begin{eqnarray} \label{pparts1} 0 \to \Omega_{C/B} \otimes_C E \to P^1_{C/B}(E) \to E \to 0 \end{eqnarray} of $E$.  The cokernel map in \eqref{pparts1} admits a natural $B$ linear section \begin{eqnarray} \label{s} s : E & \to & P^1_{C/B}(E) \\ \nonumber e & \mapsto & 1 \otimes 1 \otimes e, \end{eqnarray} which yields a $B$ module splitting of \eqref{pparts1}.  The map $s$ is generally not $C$ linear, since \begin{eqnarray*} s(c \cdot e) - c \cdot s(e) & = & 1 \otimes 1 \otimes c \cdot e - c \cdot (1 \otimes 1 \otimes e) \\ &=& 1 \otimes c \otimes e - c \otimes 1 \otimes e \\ &=& dc \otimes e. \end{eqnarray*}

Now suppose that $E = C \otimes_B M$ is obtained by extension of scalars from a $B$ module $M$.  Then I claim the map $t : E \to P^1_{C/B}(E)$ given by \begin{eqnarray} \label{t} t : C \otimes_B M & \to & ((C \otimes_B C) / I_\Delta^2) \otimes_C (C \otimes_B M) \\ \nonumber c \otimes m & \mapsto & c \otimes 1 \otimes 1 \otimes m \end{eqnarray} provides a $C$ linear splitting \begin{eqnarray} \label{splitting1} P^1_{C/B}(E) &=& (\Omega_{C/B} \otimes E) \oplus E \end{eqnarray}  of \eqref{pparts1}.  Obviously this is a splitting; it is $C$ linear by the following computation: \begin{eqnarray*} t(c' \cdot c \otimes m)- c' \cdot t(c \otimes m) &=& t(cc' \otimes m)-c' \cdot (c \otimes 1 \otimes 1 \otimes m) \\ &=& c'c \otimes 1 \otimes 1 \otimes m - c'c \otimes 1 \otimes 1 \otimes m \\ &=& 0. \end{eqnarray*} 

\subsection{} \label{redatiyah2} Continue to assume $E=C \otimes_B M$ and suppose \begin{eqnarray} \label{SES1} \xymatrix{ 0  \ar[r] & N \ar[r] & E \ar[r]^f & F \ar[r] & 0 } \end{eqnarray} is an exact sequence of $C$ modules.  From the naturality of the principal parts sequence we obtain a commutative diagram \begin{eqnarray} \label{horparts2} \xymatrix{ & 0 & 0 & 0 \\ 
0 \ar[r] & \Omega_{C/B} \otimes_C F  \ar[u] \ar[r]                          & P^1_{C/B}(F) \ar[u] \ar[r] & F \ar[u] \ar[r] & 0 \\ 
0 \ar[r] & \Omega_{C/B} \otimes_C E \ar[u]^{\Omega_{C/B} \otimes f}  \ar[r] & P^1_{C/B}(E) \ar[r] \ar[u] & E \ar[u]_f \ar[r] & 0 \\
0 \ar[r] & \Omega_{C/B} \otimes_C N \ar[u]  \ar[r]                          & P^1_{C/B}(N) \ar[r] \ar[u] & N \ar[u] \ar[r] & 0 \\ 
         &                                                                  &                            & 0 \ar[u] } \end{eqnarray} with exact rows and columns.

Using the splitting \eqref{splitting1} and the diagram \eqref{horparts2} we obtain a morphism of double complexes of $C$ modules \begin{eqnarray} \label{beta} \xymatrix{  \Omega_{C/B} \otimes_C E         \ar@{=}[r] & \Omega_{C/B} \otimes_C E \\  
            \Omega_{C/B} \otimes_C N \ar[u]  \ar[r] & P^1_{C/B}(N) \ar[u]  } & \xymatrix@C+10pt@R-10pt{ \\ \ar[r]^\beta &} & \xymatrix{  \Omega_{C/B} \otimes_C E \\  
            \Omega_{C/B} \otimes_C N \ar[u]    } \end{eqnarray} (where $\Omega_{C/B} \otimes_C N$ is placed in degree $(0,-1)$) by projecting onto the quotient complex.  The induced map on $H^0$ is given by \begin{eqnarray} \label{beta0} H^0(\beta) : N & \to & \Omega_{C/B} \otimes_C F \\ \nonumber c \otimes m & \mapsto & dc \otimes f(1 \otimes m). \end{eqnarray}  Indeed, $H^0(\beta)(c \otimes m)$ can be computed by choosing (locally) a lift of $c \otimes m$ to $P^1_{C/B}(E)$, then applying $\Omega_{C/B} \otimes f$ to the component of this lift in $\Omega_{C/B} \otimes_B E$ under the splitting \eqref{splitting1} (the result will be independent of the lift of $c \otimes m \in N \subset E$ and will yield a $C$ linear map).  We may as well choose the lift systematically be means of the $B$ linear map $s$ in \eqref{s}.  Since the splitting \eqref{splitting1} is defined via $t$, the $\Omega_{C/B} \otimes E$ component of $s(c \otimes m)$ is given by $s(c \otimes m) - t(c \otimes m)$.  Now we simply compute \begin{eqnarray*} H^0(\beta)(c \otimes m) &=& (\Omega_{C/B} \otimes f)(s(c \otimes m)-t(c \otimes m)) \\ &=&  (\Omega_{C/B} \otimes f)( 1 \otimes 1 \otimes c \otimes m - c \otimes 1 \otimes 1 \otimes m) \\ &=& (\Omega_{C/B} \otimes f)( 1 \otimes c \otimes 1 \otimes m - c \otimes 1 \otimes \otimes 1 \otimes m) \\ &=& (\Omega_{C/B} \otimes f)(dc \otimes 1 \otimes m) \\ &=& dc \otimes f(1 \otimes m). \end{eqnarray*}  Notice that this map is the same as the one in \eqref{H0ra}.
			
\subsection{}  \label{redatiyah3} In order to generalize the results of \eqref{redatiyah1} and \eqref{redatiyah2} by replacing K\"ahler differentials with the cotangent complex it becomes necessary to impose some additional hypotheses.  This is because, even though $P^1_{C/B}(E)$ splits naturally when $E$ is pulled back from $B$, there is no reason to think that $P^1_{P/B}(E)$ splits when $P=P_B C$ is the standard simplicial resolution of $C$ by free $B$ algebras and $E$ is viewed as a $P$ module by restriction of scalars along the augmentation $P \to C$.  What is needed is a particularly nice choice of resolution of $C$.  We adopt, for the remainder of this section, the setup of \eqref{setup} and we assume, as in \eqref{assumptions} that $B_1$ is flat over $A$.  To simplify notation, we will just write $N,E=B_1 \otimes_A B_2,F$ (as in the previous section) instead of $N_0,E_0,F_0$, since we will not consider here any problem of lifting module maps over algebra extensions.

Let $P_1 = P_A B_1$ be the standard simplicial resolution of $B_1$ by free $A$ algebras.  Since the homology $B_1$ of $P_1$ (viewing $P_1$ just as a complex of $A$ modules) is $A$ flat, taking homology commutes with tensoring with any $A$ module.  In particular, $P := P_1 \otimes_A B_2$ is a simplicial resolution of $B_0 = B_1 \otimes_A B_2$ by free $B_2$ algebras and hence the cotangent complex $\LL_{B_0/B_2}$ is the image in $D(B_0)$ of $\Omega_{P/B_2} \otimes_P B_0$ (II.2.1.2).  Note that the natural isomorphism $\Omega_{P/B_2} = \Omega_{P_1/A} \otimes_A B_2$ corresponds to the natural derived category isomorphism \begin{eqnarray} \label{flatiso} \LL_{B_0/B_2} &=& \LL_{B_1/A} \otimes_{B_1}^{\bL} B_0 \end{eqnarray} obtained from the base change theorem (II.2.2).  (Actually, this is more-or-less how one proves the base change properties of the cotangent complex.) Set $V = P_1 \otimes_A E_2$, so $V$ is obtained from the $B_2$ module $E_2$ by extension of scalars along $B_2 \to P$ and hence, just as in \eqref{redatiyah1}, we obtain a natural map $t : V \to P^1_{P/B_2}(V)$ yielding a splitting \begin{eqnarray} \label{splitting2} P^1_{P/B_2}(V) = \Omega_{P/B_2} \otimes_P V \oplus V \end{eqnarray} of the principal parts sequence of $V$.  The augmentation map $V \to B_1 \otimes_A E_2 = E$ is a quasi-isomorphism of $P$ modules (viewing $E$ as a $P$ module via restriction of scalars along $P \to B_0$), again because the homology of $P_1$ is flat over $A$.  We have a morphism of $P$ module extensions \begin{eqnarray} \label{extensionsplitting} & \xymatrix{ 0 \ar[r] & \Omega_{P/B_2} \otimes_P V \ar[d]_{\simeq} \ar[r] & (\Omega_{P/B_2} \otimes_P V) \oplus V \ar[d]^{\simeq} \ar[r] & V \ar[r] \ar[d]^{\simeq} & 0 \\ 0 \ar[r] & \Omega_{P/B_2} \otimes_P E \ar[r] & P^1_{P/B_2}(E) \ar[r] & E \ar[r] & 0,} \end{eqnarray} where the vertical arrows are quasi-isomorphisms (by the Five Lemma and the fact that $\Omega_{P/B_2}$ is a flat $P$ module) so that the parts sequence of $E$ is quasi-isomorphic to a split sequence of $P$ modules.

As in \eqref{redatiyah2}, the exact sequence \eqref{SES1} yields a commutative diagram \begin{eqnarray} \label{horparts4} \xymatrix{ & 0 & 0 & 0 \\ 
0 \ar[r] & \Omega_{P/B_2} \otimes_P F  \ar[u] \ar[r] & P^1_{P/B_2}(F) \ar[u] \ar[r] & F \ar[u] \ar[r] & 0 \\ 
0 \ar[r] & \Omega_{P/B_2} \otimes_P E \ar[u]  \ar[r] & P^1_{P/B_2}(E) \ar[r] \ar[u] & E \ar[u] \ar[r] & 0 \\
0 \ar[r] & \Omega_{P/B_2} \otimes_P N \ar[u]  \ar[r] & P^1_{P/B_2}(N) \ar[r] \ar[u] & N \ar[u] \ar[r] & 0 \\ 
         &         0 \ar[u]                          &      0   \ar[u]                   & 0 \ar[u] } \end{eqnarray} of $P$ modules with exact rows and columns (except here $\Omega_{P/B_2}$ is a flat $P$ module, so the complex is ``more exact").   Consider the two term complex \begin{eqnarray} W & := & [ \xymatrix{ \Omega_{P/B_2} \otimes_P F \oplus V \ar[r] & P^1_{P/B_2}(F) } ] \end{eqnarray} of $P$ modules (in degrees $0,1$), where the map $V \to P^1_{P/B_2}(F)$ is the composition of the map $V \to P^1_{P/B_2}(E)$ appearing in \eqref{extensionsplitting} and the natural map $P^1_{P/B_2}(E) \to P^1_{P/B_2}(F)$.  From the exactness of \eqref{horparts4} and the quasi-isomorphisms in \eqref{extensionsplitting} it follows that $W$ is quasi-isomorphic to $N$.\footnote{This is probably bad terminology.  We don't mean that the homology of $W$ is \emph{isomorphic} to $N$ as a $P$ module, but only that the homology of $W$ is \emph{quasi-isomorphic} to $N$ as simplicial abelian groups.  In other words, the corresponding double complex of abelian groups has homology $N$.}  This is clear once we note that $W$ is quasi-isomorphic to the double complex below.  \begin{eqnarray} \label{beta2} \xymatrix{  (\Omega_{P/B_2} \otimes_P E ) \oplus V        \ar[r] & P^1_{P/B_2}(E) \\  
            \Omega_{P/B_2} \otimes_P N \ar[u]  \ar[r] & P^1_{P/B_2}(N) \ar[u]  } \end{eqnarray}  There is an obvious morphism of complexes $\beta : W \to \Omega_{P/B_2} \otimes_P F$ obtained by projecting to the quotient complex.  After extending scalars along $P \to B_0$ and using the isomorphism \eqref{flatiso}, we may view $\beta$ as a $D(B_0)$ morphism  \begin{eqnarray} \label{reducedatiyahclassdefn} \beta: N \to \LL_{B_1/A} \otimes^{\bL}_{B_1} F  \end{eqnarray} called, suggestively, the \emph{reduced Atiyah class} of the quotient $f : E \to F$.
			
\subsection{}  We now explain the appellation ``reduced Atiyah class."  Notice that the projection map \begin{eqnarray} \label{factorization} [\xymatrix{ \Omega_{P/B_2} \otimes_P F \ar[r] & P^1_{P/B_2}(F) }] & \xymatrix{ \ar[r] & } & \Omega_{P/B_2} \otimes_P F \end{eqnarray} obviously factors through $\beta$ by including the domain complex as a subcomplex of $W$.  

By (IV.2.3.7.3) the image of this projection map in the derived category $D(B_0)$ coincides (up to a shift) with the Atiyah class \begin{eqnarray} \label{atiyahclass} {\rm at}_{B_0/B_2}(F) : F & \to & \LL_{B_0/B_2} \otimes^{\bL}_{B_0} F[1] \end{eqnarray} of $F$ relative to $B_2 \to B_0$ (see (IV.2.3) for the construction of the Atiyah class via the graded cotangent complex).  Note that the aforementioned inclusion of subcomplexes is naturally isomorphic in $D(B_0)$ to the map $F[-1] \to N$ obtained from \eqref{SES1}.  Hence, in $D(B_0)$, we obtain a commutative diagram \begin{eqnarray} \label{factorizationofatiyah} \xymatrix{ F[-1] \ar[d] \ar[rd]^-{ \at_{B_0/B_2}(F)[-1] } \\ N \ar[r]^-\beta & \LL_{B_0/B_2} \otimes^{\bL} F } \end{eqnarray}  factoring the (shift of the) Atiyah class of $F$ through the reduced Atiyah class $\beta$.

\subsection{} \label{redatiyah4} The relationship between the reduced Atiyah class of \eqref{redatiyah3} and the one introduced in \eqref{reducedatiyah} can be explained as follows.  From the transitivity triangle of graded cotangent complexes associated to the diagram of graded rings \begin{eqnarray} \label{gradedrings2} & \xymatrix{ B_2[E_2] \ar[r] & B_0[E] \ar[r] & B_0[F] \\ B_2 \ar[u] \ar[r] & B_0 \ar[r] \ar[u] & B_0[F] \ar@{=}[u] } \end{eqnarray} we obtain a commutative diagram \begin{eqnarray} \label{diagram6} & \xymatrix{ N \ar[r] & k^1 \LL_{B_0[E]/B_2[E_2]} \otimes^{\bL}_{B_0[E]} B_0[F] \\ F[-1] \ar[u] \ar[r]^-{{\rm at}_{B_0/B_2}(F)[-1]} & \LL_{B_0/B_2} \otimes_{B_0}^{\bL} F \ar[u]_{\cong} } \end{eqnarray} where the bottom horizontal arrow is the Atiyah class of $F$ relative to $B_2 \to B_0$.  Since $B_0[E] = B_2[E_2] \otimes_{B_2} B_0$ and $B_0$ is flat over $B_2$, it follows from the base change theorem (II.2.2) that the right vertical arrow is an isomorphism.

\subsection{Theorem} \emph{The diagram \begin{eqnarray} \label{diagram7} & \xymatrix{ N \ar[r] \ar[rd]_\beta & k^1 \LL_{B_0[E]/B_2[E_2]} \otimes^{\bL}_{B_0[E]} B_0[F] \\ F[-1] \ar[u] \ar[r] & \LL_{B_0/B_2} \otimes_{B_0}^{\bL} F \ar[u]_{\cong} } \end{eqnarray} obtained from \eqref{diagram6} by inserting the reduced Atiyah class is commutative.}  

\noindent \emph{Proof.} We already saw in \eqref{factorizationofatiyah} that the lower triangle commutes, so we focus now on the upper triangle.  Consistently with the previous notation, set \begin{eqnarray*} P_1 & := & P_A B_1 \\ P & := & P_1 \otimes_A B_2 \\ V &:= & P_1 \otimes_A E_2 \\ Q &:=& P_P^{\Delta {\rm gr}} P[F] \\ R &:=& P_{P[V]}^{\Delta {\rm gr}} P[F], \end{eqnarray*} where, for example, $P_P^{\Delta {\rm gr}} P[F]$ is the diagonal of the standard bisimplicial graded resolution of $P[F]$ by graded free $P$ algebras.  The diagram of graded simplicial rings \begin{eqnarray} \label{gradedrings3} & \xymatrix{ B_2[E_2] \ar[r] & P[V] \ar[r] & R \\ B_2 \ar[u] \ar[r] & P \ar[r] \ar[u] & Q \ar[u] } \end{eqnarray} admits a natural augmentation quasi-isomorphism to the diagram \eqref{gradedrings2} and the transitivity triangle of graded cotangent complexes associated to \eqref{gradedrings2} is naturally identified with the map of $D(B_0[F])$ triangles associated to the map of short exact sequences of $R$ modules \begin{eqnarray} \label{SESmap} & \xymatrix{ 0 \ar[r] & \Omega_{P[V]/B_2[E_2]} \otimes_{P[V]} R \ar[r] & \Omega_{R/B_2[E_2]} \ar[r] & \Omega_{R/P[V]} \ar[r] & 0 \\ 0 \ar[r] & \Omega_{P/B_2} \otimes_P R \ar[r] \ar@{=}[u] & \Omega_{Q/B_2} \otimes_Q R \ar[r] \ar[u] & \Omega_{Q/P} \otimes_Q R \ar[u] \ar[r] & 0 } \end{eqnarray} after extending scalars along the composition $R \to P[F] \to B_0[F]$ of the augmentation maps.  Note that the left vertical arrow in \eqref{SESmap} is an isomorphism because the left square in \eqref{gradedrings3} is cocartesian and note that the sequences are exact on the left because $R$ is free term-by-term over $P[V]$ and $Q$ is free term-by-term over $P$ (II.1.1.2.13).

The composition $N \to \LL_{B_0/B_2} \otimes^{\bL}_{B_0} F$ of the top horizontal arrow and the right vertical arrow in \eqref{diagram7} is the image in $D(B_0)$ of the quotient complex projection \begin{eqnarray} \label{firstmap}  [\xymatrix{ \Omega_{P/B_2} \otimes_P R^1 \ar[r] & k^1 \Omega_{R/B_2[E_2]} }] & \xymatrix{ & \ar[r] & } & \Omega_{P/B_2} \otimes_P R^1 \end{eqnarray} obtained from the degree one part of the diagram \eqref{SESmap} (after extending scalars from $R^0 := k^0 R$ to $B_0 = k^0 B_0[F]$).  Note that the terms in the domain complex of \eqref{firstmap} should be placed in degrees $0,1$ since we have natural isomorphisms $$k^1 \Omega_{R/P[V]} \otimes_{R^0} B_0 = k^1 \LL_{B_0[F]/B_0[E]} = N[1].$$

From the graded ring map $B_2 \to B_2[E_2] \to R$, the fact that $\Omega_{B_2[E_2]/B_2} = k^1 \Omega_{B_2[E_2]/B_2}=E_2$ and the fact that $R$ is free term-by-term over $B_2[E_2]$, we obtain an exact sequence \begin{eqnarray} \label{seqz} & \xymatrix{ 0 \ar[r] & E_2 \otimes_{B_2} R^1 \ar[r] & k^1\Omega_{R/B_2[E_2]} \ar[r] & k^1 \Omega_{R/B_2} \ar[r] & 0 } \end{eqnarray} of $R^0$ modules. From (II.1.2.6.7) and the proof of (IV.2.3.7.3) we obtain a morphism \begin{eqnarray} \label{seqx} & \xymatrix{ 0 \ar[r] & \Omega_{P/B_2} \otimes_P R^1 \ar[r] \ar@{=}[d] & P^1_{P/B_2}(R^1) \ar[r] \ar[d]_{\simeq} & R^1 \ar[d]^d_{\simeq} \ar[r] & 0 \\ 0 \ar[r] & \Omega_{P/B_2} \otimes_P R^1 \ar[r] & k^1(\Omega_{R/B_2}) \ar[r] & \Omega_{R/P} \ar[r] & 0 } \end{eqnarray} of exact sequences of $R^0$ modules. Since $P \to R$ is a quasi-isomorphism in degree zero and $P$ is supported in degree zero, the right vertical arrow is a quasi-isomorphism (IV.2.2.5), hence so is the middle vertical arrow.

We now have a sequence of quasi-isomorphisms compatible with the projections to $\Omega_{P/B_2} \otimes_P R^1$ as follows. \begin{eqnarray*} & &  [\xymatrix{ \Omega_{P/B_2} \otimes_P R^1 \ar[r] & k^1 \Omega_{R/B_2[E_2]} } ] \\ & \simeq & [ \xymatrix{ (\Omega_{P/B_2} \otimes_P R^1) \oplus (E_2 \otimes_{B_2} R^0) \ar[r] & k^1 \Omega_{R/B_2} } ] \\ & \simeq & [ \xymatrix{ (\Omega_{P/B_2} \otimes_P R^1) \oplus (E_2 \otimes_{B_2} R^0) \ar[r] & P^1_{P/B_2}(R^1) }] \\ & \simeq & [ \xymatrix{ (\Omega_{P/B_2} \otimes_P F) \oplus V \ar[r] & P^1_{P/B_2}(F) }] . \end{eqnarray*}  Note that the last quasi-isomorphism is compatible with the quasi-isomorphism $$\Omega_{P/B_2} \otimes_P R^1 \simeq \Omega_{P/B_2} \otimes_P F.$$  The first is obtained using the sequence \eqref{seqz} and the second is obtained using \eqref{seqx}.  The other quasi-isomorphisms are obtained from the natural augmentation maps: for example $R^0 \simeq P$ and $R^1 \simeq F$ are the degree $0,1$ parts of the augmentation quasi-isomorphism $R \simeq P[F]$.  This completes the proof.

\subsection{}  The reduced Atiyah class is functorial in various ways that we leave to the reader to make precise.  We only point out that, given a morphism of $B_0$ module extensions of the form \begin{eqnarray} \label{SESmorphism} \xymatrix{ 0  \ar[r] & N \ar[r] \ar[d] & E \ar@{=}[d] \ar[r]^f & F \ar[r] \ar[d] & 0 \\ 0 \ar[r] & N' \ar[r] & E \ar[r] & F' \ar[r] & 0,} \end{eqnarray} we have a commutative diagram in $D(B_0)$: \begin{eqnarray} \label{naturalityofredatiyah} \xymatrix{ N \ar[d] \ar[r]^-{\beta} & \LL_{B_0/B_2} \otimes^{\bL}_{B_0} F \ar[d] \\ N' \ar[r]^-{\beta'} & \LL_{B_0/B_2} \otimes^{\bL}_{B_0} F'} \end{eqnarray}  It is a good exercise to prove this using both constructions of the reduced Atiyah class.

\subsection{Exercise} Show that the Atiyah class factors through the reduced Atiyah class as in \eqref{atiyahclass} by using the constructions of the two classes via the graded cotangent complex (instead of using principal parts as we have done). 

\section{Perfect quotients}  \label{perfectquotients}

\subsection{}  We continue with the setup \eqref{setup} and the assumptions \eqref{assumptions}.  We also assume that the extension \eqref{extension1} is trivialized: \begin{eqnarray} \label{extension3} \underline{B}_1 & = & \xymatrix{ 0 \ar[r] & I_1 \ar[r] & B_1[I_1] \ar[r] & B_1 \ar[r] & 0 \\ & & & A, \ar[ul] \ar[u] } \end{eqnarray}  and that $F_0$ is flat over $B_1$.  

The extension \eqref{extension2} is also naturally trivialized: $B=B_0[I]$.  Given an $A$ algebra section $s : B_1 \to B_1[I_1]$ of $B_1[I_1] \to B_1$, we obtain a flat deformation (c.f.\ \eqref{flatdeformations}) $$s^* f_0 := B_1[I_1] \otimes_{s} f_0 $$ of $f_0$ over $B$.  Here we write $B_1[I_1] \otimes_{s} $ instead of $B_1[I_1] \otimes_{B_1}$ to emphasize that $B_1[I_1]$ is regarded as a $B_1$ module via restriction of scalars along $s$.  Note that $s^* E_0$ is naturally identified with $E$ via the isomorphism \begin{eqnarray} \label{naturaliso} B_1[I_1] \otimes_s (B_1 \otimes_A E_2) & \to  & B_1[I_1] \otimes_A E_2 \\ \nonumber (b_1+i_1) \otimes b_1' \otimes e_2 & \mapsto & (b_1+i_1)s(b_1') \otimes e_2  \end{eqnarray} of $B$ modules.  Set $F_s := B_1[I_1] \otimes_{s} F_0$.  The solution $s^*f_0$ to \eqref{liftingproblem} can be explicitly written \begin{eqnarray} \label{F_s} \xymatrix{ 0 \ar[r] & I_1 \otimes_{A} E_2 \ar[r] \ar[d]_{I_1 \otimes f_0} &  B_1[I_1] \otimes_A E_2 \ar[r] \ar[d] & B_1 \otimes_A E_2 \ar[d]^{f_0} \ar[r] & 0 \\ 0 \ar[r] & I_1 \otimes_{B_1} F_0 \ar[r] & F_s \ar[r] & F_0 \ar[r] & 0 } \end{eqnarray}

Since $F_0$ is flat over $B_1$, note that the kernel of $E \to F_s$ is $N_s := B_1[I_1] \otimes_s N_0$ and the kernel of $I_1 \otimes_{B_1} E_0 \to I_1 \otimes_{B_1} F_0$ is $I_1 \otimes_{B_1} N_0$.  Let $s_0 : B_1 \to B_1[I_1]$ denote the zero section $b_1 \mapsto b_1$ of $B_1[I_1] \to B_1$.

We say that the surjection $f_0 : E_0 \to F_0$ is a \emph{perfect quotient} if the map \begin{eqnarray*}  \{ A {\rm  \; algebra \; sections \; of \; } B_1[I_1] \to B_1 \; \} & \to & \{ {\rm \; flat \; deformations \; of \; } f_0 {\rm \; over \;} B_1[I_1] \; \} \\ s & \mapsto & s^* f_0 \end{eqnarray*} is bijective for every $B_1$ module $I_1$.

Given an $A$ algebra section $s$ of $B_1[I_1] \to B_1$, let $d_s : B_1 \to I_1$ be the corresponding $A$ linear derivation and let $g_s : \Omega_{B_1/A} \to I_1$ be the corresponding map of $B_1$ modules.

\subsection{Lemma} \label{mainresultlemma} \emph{ The map $N_0 \to I_1 \otimes_{B_0} F_0$ defined by the map of complexes \begin{eqnarray} \label{map1} & \xymatrix{ I_1 \otimes_{B_1} N_0 \ar[r] & N_{s} \ar[d] \\ & F_{s_0} \ar[r] & F_0. } \end{eqnarray} as in Remark~\ref{torsorbijection} coincides with the map induced by the composition \begin{eqnarray} \label{map2} & \xymatrix{ N_0 \ar[r] & \LL_{B_1/A} \otimes^{\bL} F_0 \ar[r] &  \Omega_{B_1/A} \otimes^{\bL} F_0 \ar[r] & I_1 \otimes_{B_1} F_0 } \end{eqnarray} of the truncated reduced Atiyah class of $f_0$ and $g_s \otimes F_0$.  Both maps are given by \begin{eqnarray} \label{map3} N_0 & \to & I_1 \otimes_{B_1} F_0 \\ \nonumber b_1 \otimes e_2 & \mapsto & d_s b_1 \otimes f_0(1 \otimes e_2). \end{eqnarray} }

\noindent \emph{Proof}.  The maps $s,g_s,$ and $d_s$ are related by the formulas \begin{eqnarray} s(b_1) & = & b_1+d_j(b_1) \\ \nonumber g_s(db_1) &=& d_j(b_1),\end{eqnarray} so the fact that \eqref{map2} is given by \eqref{map3} follows from \eqref{beta0}.

To calculate the map defined by \eqref{map1}, we will make use of the natural $B_2$ linear section $n \mapsto 1 \otimes n$ of $N_s \to N_0$.   Note that the isomorphism $E_s \cong E_{s_0}$ obtained by identifying both $E_s$ and $E_{s_0}$ with $E$ via the natural isomorphism \eqref{naturaliso} can be explicitly written: \begin{eqnarray} \label{naturaliso2} (B_1[I_1] \otimes_A B_2) \otimes_s (B_1 \otimes_A E_2) & \to & (B_1[I_1] \otimes_A B_2) \otimes_{s_0} (B_1 \otimes_A E_2) \\ \nonumber (b_1+i_1) \otimes b_2 \otimes b_1' \otimes e_2 & \mapsto & (b_1+i_1)s(b_1') \otimes b_2 \otimes 1 \otimes e_2. \end{eqnarray}  The map defined by \eqref{map1} can be computed on $n=b_1 \otimes e_2 \in N=\Ker f_0$ by following $1 \otimes 1 \otimes b_1 \otimes e_2$ through the sequence of maps \begin{eqnarray*} & \xymatrix{  (B_1[I_1] \otimes_A B_2) \otimes_s (B_1 \otimes_A E_2) \ar[d]^{\eqref{naturaliso2}}_{\cong} \\ (B_1[I_1] \otimes_A B_2) \otimes_{s_0} (B_1 \otimes_A E_2) \ar[d] \\  (B_1[I_1] \otimes_A E_2) \otimes_{s_0} F_0 } \end{eqnarray*} and noting that the result is in $(I_1 \otimes_A B_2) \otimes_{B_0} F_0$.  Carrying this out, we find: \begin{eqnarray*} 1 \otimes 1 \otimes b_1 \otimes e_2 & \mapsto & s(b_1) \otimes 1 \otimes 1 \otimes e_2 \\ &=& b_1 \otimes 1 \otimes 1 \otimes e_2 + d_sb_1 \otimes 1 \otimes 1 \otimes e_2 \\ &=& 1 \otimes 1 \otimes b_1 \otimes e_2 + d_s b_1 \otimes 1 \otimes 1 \otimes e_2 \\ & \mapsto & 1 \otimes 1 \otimes f_0(b_1 \otimes e_2)+d_s b_1 \otimes 1 \otimes f_0(1 \otimes e_2) \\ & = & d_s b_1 \otimes 1 \otimes f_0(1 \otimes e_2).  \end{eqnarray*} Of course we drop the $1$ in the middle when we make the natural identification $$(I_1 \otimes_A B_2) \otimes_{B_0} F_0 = I_1 \otimes_{B_1} F_0.$$  We have now shown that the map induced by \eqref{map1} is also given by \eqref{map3}, so the proof is complete.

\section{Applications to the Quot scheme} \label{applications}

\subsection{Setup}  We work throughout in the category $\Sch$ of schemes over a field $k$, which we refer to simply as ``schemes".  All relative constructions done without explicit reference to a morphism are assumed to be relative to the terminal object, so, for example, $\LL_{X} = \LL_{X/k}$, $X \times Y = X \times_k Y$, etc.\  Let $Y$ be a scheme, $E$ a quasi-coherent sheaf on $Y$.  We assume that the functor \begin{eqnarray} \label{quotfunctor} \Sch^{\rm op} & \to & \Sets \\ \nonumber X & \mapsto & \{ {\rm \; quotients \; of \; } \pi_2^* E {\rm \; on \; } X \times Y {\rm \; flat \; over \;} X \; \} \end{eqnarray} is representable by a scheme $Q$ (the \emph{Quot scheme}).  This holds, for example, if $Y$ is projective and $E$ is coherent, in which case the components of $Q$ are projective \cite{Gro}.  The Quot functor \eqref{quotfunctor} is representable by an algebraic space in much more generality.  Assuming the basic machinery of Serre duality, etc.\ for algebraic spaces, the results of this section carry over to that setting as well.

Since we work over a field, for schemes $X$ and $Y$ we have $\LL_{X \times Y} = \pi_1^* \LL_X \oplus \pi_2^* \LL_Y$ and $\LL_{X \times Y / Y} = \pi_1^* \LL_X$ (the projections are flat).  For an exact sequence \begin{eqnarray} \label{SES5} & \xymatrix{ 0 \ar[r] & N \ar[r] & \pi_2^* E \ar[r]^f & F \ar[r] & 0 } \end{eqnarray} on $X \times Y$ (with $F$ flat over $X$) we let $\beta(f) : N \to \pi_1^* \LL_X \otimes^{\bL} F$ (or just $\beta$ if $f$ is clear from context) denote the reduced Atiyah class of $f$ \eqref{reducedatiyah}.\footnote{It is not necessary to derive the tensor product in $\pi_1^* M \otimes^{\bL} F$ because $$M \mapsto \pi_1^* M \otimes F = \pi_1^{-1} M \otimes_{\pi_1^{-1} \O_X} F$$ is an exact functor from $X$ modules to $X \times Y$ modules since $F$ is flat over $X$}  Recall the natural transformation  \begin{eqnarray*} \ov{\beta} : \Hom_{D(X)}( \LL_X, \slot) & \to & \Hom_{D(X \times Y)}(N, \pi_1^*  \slot \otimes^{\bL} F) \\ \ov{\beta}(M)(g) & := & (\pi_1^*g \otimes^{\bL} F) \beta. \end{eqnarray*} from the introduction.  By abuse of notation, we write $f : X \to Q$ for the morphism corresponding to the quotient $f$.

Since all the technical results are already in place, the proof of the next theorem will amount to little more than unwinding various definitions.  I view this theorem as the main result of this paper.

\subsection{Theorem} \label{mainresult} \emph{ The following are equivalent: \begin{enumerate} \item $f$ is a perfect quotient. \item $\ov{\beta}(I)$ is an isomorphism for every quasi-coherent sheaf $I$ on $X$. \end{enumerate}  If $f:X \to Q$ is formally \'etale, then these two equivalent conditions hold, and, furthermore, $$ \ov{\beta}(I[1]) : \Ext^1(\LL_X,I) \to \Ext^1(N,\pi_1^*I \otimes F) $$ is injective for every quasi-coherent sheaf $I$ on $X$.}

\noindent \emph{Proof.}  For a quasi-coherent sheaf $I$, let $\iota_I : X \into X[I]$ be the trivial square zero closed immersion with ideal $I$.  By Lemma~\ref{mainresultlemma} we have a commutative diagram \begin{eqnarray}  & \xymatrix{ \Hom(\Omega_X,I) \ar@{=}[d] \ar[r]^-{\ov{\beta}(I)} & \Hom(N,\pi_1^*I \otimes F) \ar@{=}[d] \\  \{  {\rm  \; retracts \; of \; } \iota_I \; \} \ar[r] & \{ {\rm \; flat \; deformations \; of \; } f {\rm \; over \;} X[I] \; \}  } \end{eqnarray} where the bottom arrow is given by pullback.  By defintion, the latter is an isomorphism for every $I$ iff $f$ is a perfect quotient, so the equivalence of the two conditions is clear.

If $f:X \to Q$ is formally \'etale, then by defintion of ``formally \'etale" $s \mapsto fs$ gives a bijection between retracts of $\iota_I$ and completions of the solid diagram $$ \xymatrix{ X \ar[r]^{\iota_I} \ar[d]_f & X[I] \ar@{.>}[ld] \\ Q } $$ in $\Sch$.  On the other hand, by definition of $Q$, such completions are the same thing as flat deformations of $f$ over $X[I]$.  This proves the second statement.

For the final statement, note that the Fundamental Theorem of the Cotangent Complex (III.1.2.3) identifies $\Ext^1(\LL_X,I)$ with (isomorphism classes of) square zero thickenings $X \into X'$ of $X$ with ideal sheaf $I$.  According to \eqref{atiyahclassthm}, the map $$\ov{\beta}(I[1]) : \Ext^1(\LL_X,I) \to \Ext^1(N,\pi_1^*I \otimes F)$$ takes the class $e(X') : \LL_X \to I[1]$ of the thickening $X \into X'$ to the obstruction  $$ \omega = (\pi_1^* e(X') \otimes F) \beta $$ to finding a flat deformation of $f$ over $X'$.  But by definition of the Quot scheme $Q$, finding such a flat deformation is the same thing as finding a commutative diagram of solid arrows $$\xymatrix{ X \ar@{=}[r] \ar@{^(->}[d] & X \ar[d]^f \\ X' \ar[r] \ar@{.>}[ru] & Q }$$ and by defintion of formally \'etale every such diagram can be completed to a commutative diagram as indicated by the dotted arrow.  Such a dotted arrow is the same thing as a retract of $X \into X'$, which is the same thing as an identification of $X'$ with the trivial thickening: $X'=X[I]$.  If $ \ov{\beta}(I[1])$ failed to be injective we would therefore find some nontrivial thickening $X'$ with a retract, which is absurd.

\subsection{Lemma}  \label{lemma2} \emph{ Let $Z$ be a scheme.  Suppose $F \in D(Z)$ is a complex with finite perfect amplitude whose rank $r$ is prime to the characteristic of $k$.  Let $L,N$ be arbitrary complexes of $\O_X$ modules and let $\beta : N \to L \otimes^{\bL} F$ be a morphism in $D(Z)$.  Set $F^\lor := \bR \sHom(F,\O_Z)$.  For any $E \in D(Z)$, there is an isomorphism $$\Hom_{D(Z)}(N,E \otimes^{\bL} F) \to \Hom_{D(Z)}(N \otimes^{\bL} F^\lor, E),$$ natural in $E$. }

\noindent \emph{Proof.} Since $F$ has finite perfect amplitude, we have $F \otimes^{\bL} F^\lor = \bR \sHom(F,F)$, and we have a trace morphism $${\rm tr} : \bR \sHom(F,F) \to \O_Z$$ and a ``scalar multiplication" morphism $${\rm id} : \O_Z \to \bR \sHom(F,F)$$ satisfying $({\rm tr})({\rm id}) = \cdot r$.  The desired isomorphism is given by $g \mapsto (E \otimes {\rm tr})(g \otimes F^\lor)$.  Its inverse is given by $h \mapsto (h \otimes F)(N \otimes r^{-1} {\rm id})$.

\subsection{Lemma} \label{lemma3} \emph{ Let $\Ab$ be an abelian category with enough injectives and let $\B$ be a full abelian subcategory of $\Ab$.  Let $f: E \to F$ be a map in $D(\Ab)$ between complexes with cohomology in $\B$ and vanishing in positive degrees.  Then the following are equivalent: \begin{enumerate} \item $H^0(f)$ is an isomorphism. \item $\Hom(F,I) \to \Hom(E,I)$ is an isomorphism for every $I$ in $\B$. \end{enumerate}  If these equivalent conditions are satisfied, then  the following are equivalent: \begin{enumerate} \item $H^{-1}(f)$ is surjective. \item $\Ext^1(F,I) \to \Ext^1(E,I)$ is injective for every $I$ in $\B$. \end{enumerate}}

\noindent \emph{Proof.}  This is standard.  Let $I \to J$ be an injective resolution of $I$.  Then we can compute $\Ext^n(E,I)$ as $H^n$ of the double complex $\Hom^\bullet(E,J)$.  By first taking cohomology in the $E$ direction, we obtain a spectral sequence $$\Ext^p(H^{-q}(E),I) \Rightarrow \Ext^{p+q}(E,I)$$ natural in $E,I$.  Since $H^{>0}(E)=H^{>0}(F)=0$, the map $$f^* : \Hom(F,I) \to \Hom(E,I)$$ is identified with the $\B$ morphism $$H^0(f)^* : \Hom(H^0(F),I) \to \Hom(H^0(E),I),$$ so the first statement follows from Yoneda's Lemma.  From the naturality of the low order terms in this spectral sequence, we obtain a commutative diagram \begin{eqnarray*} & \xymatrix{ 0 \ar[r] & \Ext^1(H^0(E),I) \ar[r] & \Ext^1(E,I) \ar[r] & \Hom(H^{-1}(E),I) \ar[r] & \Ext^2(H^0(E),I)  \\ 0 \ar[r] & \Ext^1(H^0(F),I) \ar[u]^{H^0(f)^*}_{\cong} \ar[r] & \Ext^1(F,I) \ar[u]^{f^*} \ar[r] & \Hom(H^{-1}(F),I) \ar[u]^{H^{-1}(f)^*} \ar[r] & \Ext^2(H^0(F),I) \ar[u]^{H^0(f)^*}_{\cong}  } \end{eqnarray*} with exact rows and natural in $I$.  Assuming $H^0(f)$ is an isomorphism, so are the indicated arrows.  By the Subtle Five Lemma,\footnote{This isn't so subtle since we have much more exactness than necessary for the desired conclusion.} $H^{-1}(f)^*$ is injective iff $f^*$ is injective.  On the other hand, $H^{-1}(f)^*$ is injective for every $I$ in $\B$ iff the $\B$ morphism $H^{-1}(f)$ is surjective.

\subsection{Remark}  One could also use the long exact sequence obtained by applying $\Hom( \slot, I)$ to the map of truncation triangles $$\xymatrix{ \tau_{ \leq -1}E[1] \ar[d] \ar[r] & E \ar[d] \ar[r] & \tau_{\geq 0} E \ar[d] \\  \tau_{ \leq -1}F[1]  \ar[r] & F  \ar[r] & \tau_{\geq 0} F  }$$ instead of the low order terms in the spectral sequence, though this amounts to the same thing.

\subsection{Theorem}  \label{secondresult} \emph{Suppose $Y$ is a projective Gorenstein scheme of pure dimension $d$, $E$ is a coherent sheaf on $Y$, and $f : \pi_2^*E \to F$ is a surjection of sheaves on $X \times Y$ with $F$ flat over $X$ and of finite perfect amplitude.  Assume that the rank of $F$ is prime to the characteristic of $k$ and that $N := \Ker f$ is also of finite perfect amplitude.  Then the functor \begin{eqnarray*} D(X) & \to & \Vect_k \\ I & \mapsto & \Hom_{D(X \times Y)}(N, \pi_1^* I \otimes^{\bL} F) \end{eqnarray*} is represented by $$\EE :=  \bR \sHom( \bR \pi_{1*} \bR \sHom(N,F) , \O_X) \in D(X).$$  If $f$ defines a formally \'etale map from $X$ to the Quot scheme, then the map $H^0(\EE \to \LL_X)$ is an isomorphism and the map $H^{-1}(\EE \to \LL_X)$ is surjective.}

\noindent \emph{Proof.}  By the hypotheses on $Y$, there is an invertible sheaf $\omega_Y$ on $Y$ and an isomorphism \begin{eqnarray} \label{Serreduality} \bR \sHom( \bR \pi_{1*} A, B) = \bR \pi_{1*} \bR \sHom(A,  \pi_1^* B \otimes \pi_2^* \omega_Y [d] ) \end{eqnarray} (Grothendieck-Serre duality) in $D(X)$ natural in $A \in D(X \times Y)$ and $B \in D(X)$.  There are natural isomorphisms \begin{eqnarray*} \Hom_{D(X \times Y)}(N, \pi_1^*I \otimes^{\bL} F) & = & \Hom_{D(X \times Y)}(N \otimes^{\bL} F^\lor, \pi_1^* I) \\ &=& \Hom_{D(X \times Y)}(N \otimes^{\bL} F^\lor \otimes \pi_2^* \omega_Y[d], \pi_1^* I \otimes \pi_2^* \omega_Y[d]) \\ &=& \Hom_{D(X)}(\bR \pi_{1*}( N \otimes^{\bL} F^\lor \otimes \pi_2^* \omega_Y[d]), I), \end{eqnarray*} where the first isomorphism is obtained from Lemma~\ref{lemma2}, the second isomorphism simply reflects the fact that $\otimes \pi_2^* \omega_Y[d]$ is an automorphism of $D(X \times Y)$ since $\pi_2^* \omega_Y$ is an invertible sheaf, and the third isomorphism is obtained from the Serre duality isomorphism \eqref{Serreduality} by applying $\bR \Gamma$ and taking $H^0$. 

Using Serre duality and the various perfection hypotheses, we obtain a sequence of natural isomorphisms \begin{eqnarray*} \EE &=&  \bR \sHom( \bR \pi_{1*} \bR \sHom(N,F) , \O_X) \\ &=& \bR \pi_{1*} ( \bR \sHom( \bR \sHom(N,F), \omega_Y[d])) \\ &=& \bR \pi_{1*}( \bR \sHom( \bR \sHom(N,F), \O_{X \times Y}) \otimes \pi_2^* \omega_Y[d]) \\ &=& \bR \pi_{1*}( \bR \sHom(N^\lor \otimes^{\bL} F, \O_{X \times Y}) \otimes \pi_2^* \omega_Y[d]) \\ &=& \bR \pi_{1*}( N \otimes^{\bL} F^\lor \otimes \pi_2^* \omega_Y[d]). \end{eqnarray*}  Putting this together with the isomorphism from the previous paragraph proves the first part of the theorem.

The second statement follows from Theorem~\ref{mainresult} and Lemma~\ref{lemma3} (applied with $\Ab = \Mod_Y$, $\B =$ quasi-coherent sheaves).  

\subsection{Lemma}  \emph{Suppose $Y=C$ is a projective curve.  Then for any quasi-compact scheme $X$ and any exact sequence \eqref{SES5} with $F$ flat on $X \times C$ and $N$ locally free, the complex $$\EE = \bR \sHom( \bR \pi_{1*} \bR \sHom(N,F) , \O_X) \in D(X)$$ is of perfect amplitude $\subseteq [-1,0]$.}

\noindent \emph{Proof.}  It suffices to show that the complex $ \bR \pi_{1*} \bR \sHom(N,F)$ is of perfect amplitude $\subseteq [0,1]$.  Since $N$ is locally free we have $\bR \sHom(N,F) = \sHom(N,F)$.  Let $D \subset C$ be an effective Cartier divisor.  Then we have an exact sequence \begin{eqnarray} & \label{SES6} \xymatrix@C-8pt{ 0 \ar[r] & \sHom(N,F) \ar[r] & \sHom(N,F)(D) \ar[r] & \sHom(N,F)(D)|_{X \times D} \ar[r] & 0 } \end{eqnarray} on $X \times C$ (writing $D$ for $\pi_2^*D$ to save notation).  Since $F$ is flat over $X$ and $X$ is quasi-compact, it follows from the basic cohomology and base change theorems that for any sufficiently positive such $D$ we have $R^{>0} \pi_{1*} \sHom(N,F)(D) = 0$.  In this case, by Grauert's Criterion, $\pi_{1*} \sHom(N,F)(D)$ will be a vector bundle on $X$.  Since $\sHom(N,F)(D)|_{X \times D}$ is supported in relative dimension zero over $X$, we conclude similarly that $$R^{>0} \pi_{1*} \sHom(N,F)(D)|_{X \times D} = 0$$ and that $\pi_{1*} \sHom(N,F)(D)|_{X \times D}$ is a vector bundle on $X$.  From the usual spectral sequences argument we conclude that $$\bR \pi_{1*} \bR \sHom(N,F) = \pi_{1*} \sHom(N,F)(D) \to \pi_{1*} \sHom(N,F)(D)|_{X \times D} $$ in $D(X)$, where the two term complex on the right is placed in degrees $0,1$.  Since the two terms in this complex are vector bundles, the proof is complete.

\subsection{Corollary}  \emph{Suppose $C$ is a projective Gorenstein curve, $E$ is a vector bundle on $C$, and $X$ is an open subset of Quot $E$ on which the universal kernel $N \subseteq \pi_2^* E$ is locally free (this always holds if $C$ is smooth).  Then the map $$\bR \sHom( \bR \pi_{1*} \bR \sHom(N,F) , \O_X) \to \LL_X$$ obtained from the reduced Atiyah class defines a perfect obstruction theory on $X$.}

\noindent \emph{Proof.}  Note that if $C$ is a smooth curve, $\O_C$ is a sheaf of PIDs so any (quasi-coherent) subsheaf of a vector bundle is a vector bundle.  If $N$ is locally free, then from the exact sequence \eqref{SES5} it is clear that $F$ has finite perfect amplitude ($\subseteq [-1,0]$ in fact), so this follows immediately from the previous results of this section.

\subsection{Remark}  If $Y$ is smooth, then $\pi_1 : X \times Y \to X$ is also smooth, so flatness of $F$ (hence $N$) over $X$ implies perfection of $F$ and $N$ by [SGA6.III.3.6].  It may be possible to weaken some of the perfection hypotheses.

\end{document}